\documentclass{icmart}
\newtheorem{thm}{Theorem}[section]

\newtheorem{cor}[thm]{Corollary}
\newtheorem{prop}[thm]{Proposition}
\newtheorem{defi}[thm]{Definition}
\newtheorem{conj}[thm]{Conjecture}
\newtheorem{question}[thm]{Question}
\theoremstyle{definition}
\newtheorem{rem}[thm]{Remark}
\newtheorem{example}[thm]{Example}

\input xypic
\xyoption{all}

\def\hbar{\bar{h}}

\def\iso{\buildrel \sim\over\to}

\def\GS{{\mathfrak{S}}}

\def\Gg{{\mathfrak{g}}}

\def\Ggl{{\mathfrak{gl}}}

\def\Gsl{{\mathfrak{sl}}}
\def\sl{{\mathfrak{sl}}}

\def\CA{{\mathcal{A}}}
\def\CB{{\mathcal{B}}}
\def\CC{{\mathcal{C}}}

\def\CE{{\mathcal{E}}}
\def\CF{{\mathcal{F}}}

\def\CI{{\mathcal{I}}}

\def\CO{{\mathcal{O}}}
\def\CP{{\mathcal{P}}}

\def\CS{{\mathcal{S}}}
\def\CT{{\mathcal{T}}}

\def\BA{{\mathbf{A}}}
\def\BB{{\mathbf{B}}}
\def\BC{{\mathbf{C}}}

\def\BF{{\mathbf{F}}}
\def\BG{{\mathbf{G}}}

\def\BL{{\mathbf{L}}}

\def\BP{{\mathbf{P}}}
\def\BQ{{\mathbf{Q}}}
\def\BR{{\mathbf{R}}}

\def\BT{{\mathbf{T}}}
\def\BU{{\mathbf{U}}}

\def\BX{{\mathbf{X}}}

\def\BZ{{\mathbf{Z}}}

\def\Br{{\mathbf{r}}}

\def\Bw{{\mathbf{w}}}

\def\Bpi{{\mathbf{\pi}}}

\def\Aut{\operatorname{Aut}\nolimits}

\def\Br{\operatorname{Br}\nolimits}

\def\cone{\operatorname{cone}\nolimits}

\def\DPic{\operatorname{DPic}\nolimits}

\def\End{\operatorname{End}\nolimits}

\def\GL{\operatorname{GL}\nolimits}
\def\gldim{\operatorname{gldim}\nolimits}

\def\Hilb{\operatorname{Hilb}\nolimits}

\def\Hom{\operatorname{Hom}\nolimits}

\def\Id{\operatorname{Id}\nolimits}

\def\Ind{\operatorname{Ind}\nolimits}

\def\Inn{\operatorname{Inn}\nolimits}

\def\Irr{\operatorname{Irr}\nolimits}

\def\ll{\operatorname{Loewy\ length}\nolimits}

\def\mMod{\operatorname{\!-mod}\nolimits}

\def\mMOD{\operatorname{\!-Mod}\nolimits}

\def\mModgr{\operatorname{\!-modgr}\nolimits}

\def\opp{{\operatorname{opp}\nolimits}}

\def\mperf{\operatorname{\!-perf}\nolimits}

\def\mlperm{\operatorname{\!-lperm}\nolimits}

\def\Pic{\operatorname{Pic}\nolimits}

\def\mproj{\operatorname{\!-proj}\nolimits}

\def\Out{\operatorname{Out}\nolimits}

\def\rank{\operatorname{rank}\nolimits}

\def\repdim{\operatorname{repdim}\nolimits}

\def\Res{\operatorname{Res}\nolimits}

\def\SL{\operatorname{SL}\nolimits}

\def\Spec{\operatorname{Spec}\nolimits}

\def\mstab{\operatorname{\!-\overline{mod}}\nolimits}

\def\StPic{\operatorname{StPic}\nolimits}

\def\TrPic{\operatorname{TrPic}\nolimits}

\def\ie{{\em i.e.}}

\def\tt{{\tilde{t}}}

\def\tR{{\tilde{R}}}

\def\idun{\mathbf{1}}
\title{Derived equivalences and finite dimensional algebras}
\author{Rapha\"el Rouquier}
\contact[rouquier@maths.leeds.ac.uk]{Rapha\"el Rouquier, Department of Pure
Mathematics, University of Leeds, Leeds, LS2 9JT, UK, and
Institut de Math\'ematiques de Jussieu, 2 place Jussieu, 75005 Paris, France.}

\date\today
\begin{document}

\begin{abstract}
We discuss the homological algebra of
representation theory of finite dimensional algebras and finite groups. We
present various methods for the construction and the study of equivalences
of derived categories: local group theory, geometry and categorifications.
\end{abstract}
\maketitle

\maketitle
\setcounter{tocdepth}{3}
\tableofcontents
\section{Introduction}

This paper discusses derived equivalences, their construction and their use,
for finite dimensional algebras, with a special focus on finite group
algebras.

\smallskip
In a first part, we discuss Brou\'e's abelian defect group conjecture and
its ramifications. This is one of the deepest problem in the representation
theory of finite groups. It is part of local representation theory, which
aims to relate
characteristic $p$ representations of a finite group with
representations of local subgroups (normalizers of non-trivial $p$-subgroups).
We have taken a more functorial viewpoint in the definition of
classical concepts (defect groups, subpairs,...).

In \S \ref{secconj},
we present Alperin's conjecture, which gives a prediction for
the number of simple representations, and Brou\'e's conjecture, which is
a much more precise prediction for the derived category, but does apply
only to certain blocks (those with abelian defect groups).

We discuss in \S \ref{secvarious} various types of equivalences that arise and
present
the crucial problem of lifting stable equivalences to derived equivalences.

In \S \ref{secloc}, we present some local methods.
We give a stronger version of the abelian
defect group conjecture that can be approached inductively and
reduced to the problem explained above of lifting stable equivalences
to derived equivalences. Roughly speaking, in a minimal counterexample
to that refinement of the abelian defect conjecture, there is a
stable equivalence.
Work of Rickard suggested to impose conditions on the terms of the complexes:
they should be direct summands of permutation modules. We explain that
one needs also to put conditions on the maps, that make the complexes look like
complexes of chains of simplicial complexes.

There is no understanding on how to construct candidates complexes who would
provide the derived equivalences expected by the abelian defect group
conjecture in general. For finite groups of Lie type (in non-describing
characteristic), we explain (\S \ref{secChev})
Brou\'e's idea that such complexes should
arise as complexes of cohomology of Deligne-Lusztig varieties. We describe
(\S \ref{secJordan})
the Jordan decomposition of blocks (joint work with Bonnaf\'e), as conjectured
by Brou\'e: Morita equivalences between blocks are constructed from
the cohomology of Deligne-Lusztig varieties. For $\GL_n$, every block is
shown to be Morita equivalent to a unipotent block. This provides some
counterpart to the Jordan decomposition of characters (Lusztig). In
\S \ref{secdlADC} and \ref{secreg}, we explain the construction of complexes
in the setting of the abelian defect conjecture. There are some delicate
issues related to the choice of the Deligne-Lusztig variety and the extension
of the action of the centralizer of a defect group to that of the normalizer.
This brings braid groups and Hecke algebras of complex reflection groups.

In \S \ref{secgeo}, we explain how to view the problem of lifting
stable equivalences to derived equivalences
as a non-commutative version of the birational
invariance of derived categories of Calabi-Yau varieties.

In \S \ref{secperv}, we describe a class of derived equivalences which
are filtered shifted Morita equivalences (joint work with Chuang). We believe
these are the building bricks for most equivalences and the associated
combinatorics should be interesting.

\smallskip

Part \S \ref{secinvar} is devoted to some invariants of
derived equivalences. In \S \ref{secauto}, we explain a functorial approach
to outer automorphism groups of finite dimensional algebras and deduce that
their identity component is preserved under various equivalences. This
functorial
approach is similar to that of the Picard group of smooth projective schemes
and we obtain also an invariance of the identity component of the product
of the Picard group by the automorphism group, under derived equivalence.

In \S \ref{secgrad}, we explain how to transfer gradings through derived
or stable equivalences. As a consequence, there should be very interesting
gradings on blocks with abelian defect. This applies as well to Hecke
algebras of type $A$ in characteristic $0$, where we obtain gradings which
should be related to geometrical gradings.

Finally, in \S \ref{secdim}, we explain the notion of dimension for
triangulated categories, in particular for derived categories of
algebras and schemes. This applies to answer a question of Auslander on the
representation dimension and a question of Benson on Loewy length of group
algebras.

\smallskip
Part \S \ref{seccat} is devoted to ``categorifications''. Such ideas have
been advocated by I.~Frenkel and have already shown their relevance in the
work of Khovanov \cite{Kh2} on knot invariants. Our idea is that 
``classical'' structures have natural
higher counterparts. These act as symmetries
of categories of representations or of sheaves.

In \S \ref{secsl2},
we explain the construction with Chuang of a categorification of $\Gsl_2$ and
we develop the associated ``$2$-representation theory''. There is an action
on the sum of
module categories of symmetric groups, and we deduce the existence of
derived equivalences between blocks with isomorphic defect groups, using
the general theory that provides a categorification of the adjoint action of
the Weyl group. This applies as well to general linear groups, and gives
a solution to the abelian defect group conjecture for symmetric and general
linear groups.

In \S \ref{secbraid}, we define categorifications of braid groups. This is
based on Soergel's bimodules.

\smallskip
I thank C\'edric Bonnaf\'e and Hyohe Miyachi for useful comments on a
preliminary version of this paper.

\section{Brou\'e's abelian defect group conjecture}

\subsection{Introduction}
\label{secintro}
\subsubsection{Blocks}
Let $\ell$ be a prime number. Let $\CO$ be the ring of integers of a finite
extension $K$ of the field $\BQ_\ell$ of $\ell$-adic numbers and $k$ its
residue field.

Let $G$ be a finite group. Modular representation theory is the
study of the categories $\CO G\mMod$ and $kG\mMod$ (finitely generated
modules).
The decomposition of $\Spec Z(\CO G)$ into connected components corresponds
to the decomposition $Z(\CO G)=\prod_b Z(\CO G)b$, where $b$ runs
over the set of
primitive idempotents of $Z(\CO G)$ (the {\em block idempotents}).
We have corresponding decompositions in {\em blocks}
$\CO G=\prod_{b} \CO G b$ and $\CO G\mMod=\bigoplus_b \CO Gb\mMod$.

\begin{rem}
One assumes usually that $K$ is big enough so that $KG$ is a product of
matrix algebras over $K$ (this will be the case if $K$ contains the
$e$-th roots of unity, where $e$ is the exponent of $G$). Descent
methods often allow a reduction to that case.
\end{rem}

\subsubsection{Defect groups}
A {\em defect group} of a block $\CO Gb$ is a minimal subgroup $D$ of $G$
such that $\Res_D^G=\CO Gb\otimes_{\CO Gb}-:
D^b(\CO Gb)\to D^b(\CO D)$ is faithful (\ie, injective on
$\Hom$'s). Such a subgroup is an $\ell$-subgroup and it is unique up
to $G$-conjugacy.

The {\em principal block} $\CO G b_0$ is the one through which the trivial
representation factors. Its defect groups are the Sylow $\ell$-subgroups of
$G$.

Defect groups measure the representation type of the block:
\begin{itemize}
\item $kGb$ is simple if and only if $D=1$
\item $kGb\mMod$ has finitely many indecomposable objects (up to isomorphism)
if and only if the defect groups are cyclic
\item $kGb$ is tame (\ie, indecomposable modules are classifiable in a
reasonable sense) if and only if the defect groups are cyclic or $\ell=2$ and
defect groups are dihedral, semi-dihedral or generalized quaternion groups.
\end{itemize}

\subsubsection{Brauer correspondence}
\label{secBrauer}
Let $\CO Gb$ be a block and $D$ a defect group. There is a unique
block idempotent $c$ of $\CO N_G(D)$ such that the restriction functor
$\Res_D^G=c\CO Gb\otimes_{\CO Gb}-:D^b(\CO Gb)\to D^b(\CO N_G(D)c)$ is faithful.

This correspondence provides a bijection
between blocks of $\CO G$
with defect group $D$ and blocks of $\CO N_G(D)$ with defect group $D$.

\subsubsection{Conjectures}
\label{secconj}
We have seen in \S \ref{secBrauer}
that $D^b(\CO G)$ embeds in $D^b(\CO N_G(D)c)$. The {\em abelian defect
conjecture} asserts that, when $D$ is abelian, the categories are actually
equivalent (via a different functor):

\begin{conj}[Brou\'e]
\label{ADC}
If $D$ is abelian, there is an equivalence $D^b(\CO Gb)\iso D^b(\CO N_G(D)c)$.
\end{conj}

A consequence of the conjecture is an isometry
$K_0(KGb)\iso K_0(KN_G(D)c)$ with good arithmetical properties (a 
{\em perfect isometry}). Note that the conjecture also carries homological
information: if $\CO Gb$ is the principal block and the equivalence sends
the trivial module to the trivial module, we deduce that the
cohomology rings of $G$ and $N_G(D)$ are isomorphic, a classical and easy fact.
It is unclear whether there should be some canonical equivalence in
Conjecture \ref{ADC}.

\smallskip
Local representation theory is the study of the relation between
modular representations and local structure of $G$. Alperin's conjecture
asserts that the number of simple modules in a block can be computed in terms
of local structure.

\begin{conj}[Alperin]
Assume $D\not=1$. Then,
$$\rank K_0(kGb)=\sum_\CS (-1)^{l(\CS)+1} \rank K_0(kN_G(\CS)c_\CS)$$
where $\CS$ runs over the conjugacy classes of chains of subgroups
$1<Q_1<Q_2<\cdots<Q_n\le_G D$, $l(\CS)=n\ge 1$ and $c_\CS$ is the sum of the
block idempotents of $N_G(\CS)$ corresponding to $b$.
\end{conj}

\begin{rem}
We have stated here Kn\"orr-Robinson's reformulation of the conjecture
\cite{KR}.
Note that the conjecture is expected to be compatible with
$\ell$-local properties of character degrees, equivariance, rationality
(Dade, Robinson, Isaacs, Navarro).
When $D$ is abelian, Alperin's conjecture (and its refinements) follows
immediately from Brou\'e's conjecture.
It would be extremely interesting to find
a common refinement of Alperin and Brou\'e's conjectures. For principal blocks,
it should contain the description of the cohomology ring as stable
elements in the cohomology ring of a Sylow subgroup.
\end{rem}

\subsection{Various equivalences}
\label{secvarious}
Let $A$ and $B$ be two symmetric algebras over a noetherian
commutative ring $\CO$.

\subsubsection{Definitions}
Let $M$ be a bounded complex of finitely generated $(A,B)$-bimodules which
are projective as $A$-modules and as right $B$-modules. Assume
there is an $(A,A)$-bimodule $R$ and a $(B,B)$-bimodule $S$ with
\begin{align*}
M\otimes_B M^*&\simeq A\oplus R\text{ as complexes of }
(A,A)\!-\!\text{bimodules}\\
M^*\otimes_A M&\simeq B\oplus S\text{ as complexes of }
(B,B)\!-\!\text{bimodules}.
\end{align*}

We say that $M$ induces a
\begin{itemize}
\item {\em Morita equivalence} if $M$ is concentrated in degree $0$ and $R=S=0$
\item {\em Rickard equivalence} if $R$ and $S$ are homotopy equivalent to
$0$ as complexes of bimodules
\item {\em derived equivalence} if $R$ and $S$ are acyclic
\item {\em stable equivalence} (of Morita type) if $R$ and $S$ are homotopy
equivalent to bounded complexes of projective bimodules.
\end{itemize}

Note that Morita $\Rightarrow$ Rickard $\Rightarrow$ stable and
Rickard $\Rightarrow$ derived. Note also
that if there is a complex inducing a stable equivalence, then there is
a bimodule inducing a stable equivalence. Finally, Rickard's theory says
that if there is a complex inducing a derived equivalence, then there is a
complex inducing a Rickard equivalence.

\smallskip
The definitions amount to requiring that $M\otimes_B -$ induces an
equivalence
\begin{itemize}
\item (Morita) $B\mMod\iso A\mMod$
\item (Rickard) $K^b(B\mMod)\iso K^b(A\mMod)$
\item (derived) $D^b(B)\iso D^b(A)$
\item (stable) $B\mstab\iso A\mstab$ (assuming $\CO$ regular)
\end{itemize}
where $K^b(A\mMod)$ is the homotopy category of bounded complexes of objects
of $A\mMod$ and $A\mstab$ is the stable category, additive quotient of
$A\mMod$ by modules of the form $A\otimes_\CO V$ with $V\in \CO\mMod$
(it is equivalent to $D^b(A)/A\mperf$ when $\CO$ is regular).

\subsubsection{Stable equivalences}
Stable equivalences arise fairly often in modular representation theory.
For example,
assume the Sylow $\ell$-subgroups of $G$ are TI, \ie, given $P$ a Sylow
$\ell$-subgroup, then $P\cap gPg^{-1}=1$ for all $g\in G-N_G(P)$.
Then, $M=\CO G$ induces a stable equivalence between $\CO G$ and
$\CO N_G(P)$, the corresponding functor is restriction (this is an
immediate application of Mackey's formula).
This restricts to a stable equivalence
between principal blocks. Unfortunately, we don't know how to derive
much numerical information from a stable equivalence.

\smallskip
A classical outstanding conjecture in representation theory of finite
dimensional algebras is
\begin{conj}[Alperin-Auslander]
\label{stablesimple}
Assume $\CO$ is an algebraically closed field.
If $A$ and $B$ are stably equivalent, then they
have the same number of isomorphism classes of simple non-projective modules.
\end{conj}

A very strong generalization of Conjecture \ref{stablesimple} is
\begin{question}
\label{liftability}
Let $A$ and $B$ be blocks with abelian defect groups and $M$ a complex of
$(A,B)$-bimodules
inducing a stable equivalence. Assume $K$ is big enough. Does there exist
$\tilde{M}$ a complex of $(A,B)$-bimodules
inducing a Rickard equivalence and such that $M$ and $\tilde{M}$
are isomorphic in $(A\otimes B^\opp)\mstab$?
\end{question}

As will be explained in \S \ref{secgluing}, this is the key step for an
inductive approach to Brou\'e's conjecture.

\begin{rem}
There are examples of blocks with non abelian defect for
which Question \ref{liftability} has a negative answer, for example
$A$ the principal block of $\mathrm{Suz}(8)$, $\ell=2$, and $B$ the principal
block of the normalizer of a Sylow $2$-subgroup (TI case), cf
\cite[\S 6]{BrLuminy}. A major problem with Question \ref{liftability}
and with Conjecture \ref{ADC} is to understand the relevance of the assumption
that the defect groups are abelian. Cf \S \ref{secgradabelian}
for a possible idea.
\end{rem}

\subsection{Local theory}
\label{secloc}
In an ideal situation, equivalences would arise from permutation modules
or more generally, from chain complexes of simplicial complexes $X$ acted on by
the groups under consideration. Then, taking fixed points on $X$
by an $\ell$-subgroup $Q$
would give rise to equivalences between blocks of the centralizers of $Q$.
We would then have a compatible system of equivalences, corresponding to
subgroups of the defect group. At the level of characters, Brou\'e defined
a corresponding notion of ``isotypie'' \cite{BrLuminy}: values of characters
at $\ell$-singular elements are related.

\subsubsection{Subpairs}
We explain here some classical facts.

A $kG$-module of the form $k\Omega$ where $\Omega$ is a $G$-set is
a permutation module.
An $\ell$-permutation module is a direct summand of a permutation module and
we denote by $kG\mlperm$ the corresponding full subcategory of $kG\mMod$.

Let $Q$ be an $\ell$-subgroup of $G$. We define the functor
$\Br_Q:kG\mlperm\to k(N_G(Q)/Q)\mlperm$:
 $\Br_Q(M)$ is the image of $M^Q$ in
$M_Q=M/\sum_{x\in Q}(x-1)M$. If $M=k\Omega$, then
$k(\Omega^Q)\iso \Br_Q(M)$:
the Brauer construction extends the fixed point construction on sets to
$\ell$-permutation modules. Note that this works only because $Q$ is an
$\ell$-group and $k$ has characteristic $\ell$.

\smallskip
To deal with non principal blocks,
we need to use Alperin-Brou\'e's subpairs.
A subpair of $G$ is a pair $(Q,e)$, where $Q$ is an $\ell$-subgroup of $G$ and
$e$ a block idempotent of $kC_G(Q)$. If we restrict to the case where $e$
is a principal block, we recover theory of $\ell$-subgroups of $G$.

A maximal subpair is of the form $(D,b_D)$, where
$D$ is a defect group of a block $kGb$ and $b_D$ is a block idempotent
of $kC_G(D)$ such that $b_Dc\not=0$ (we say that $(D,b_D)$ is a $b$-subpair).
Fix such a maximal subpair.
The $(kG,kN_G(D,b_D))$-bimodule $bkGb_D$ has, up to isomorphism,
a unique indecomposable direct summand $X$ with $\Br_{\Delta D}(X)\not=0$.
Here, we put $\Delta D=\{(x,x^{-1})\}_{x\in D}\le D\times D^\opp$.
More generally, given $\phi:Q\to R$, we put
$\Delta_\phi(Q)=\{(x,\phi(x)^{-1})\}_{x\in Q}\le Q\times R^\opp$.

We define the Brauer category $\CB r(D,b_D)$: its objects are
subpairs $(Q,b_Q)$ with $Q\le D$ and $b_Q \Br_{\Delta Q}(X)\not=0$, and
$\Hom((Q,b_Q),(R,b_R))$ is the set of $f\in\Hom(Q,R)$ such that there
is $g\in G$ with $(Q^g,b_Q^g)\in \CB r(D,b_D)$ and
 $f(x)=g^{-1}xg$ for all $x\in Q$.

\smallskip
Let $M\in kG\mlperm$ indecomposable.
A vertex-subpair of $M$ is a
subpair $(Q,b_Q)$ maximal such that $b_Q\Br_Q(M)\not=0$ (such a subpair is
unique up to conjugacy).

\subsubsection{Splendid equivalences}
Let $G$ and $H$ be two finite groups and $b$ and $b'$ two block
idempotents of $kG$ and $kH$.

The following Theorem \cite{charl,local} shows that a stable equivalence
corresponds
to ``local'' Rickard equivalences, for complexes of $\ell$-permutation modules.

\begin{thm}
\label{splendid}
Let $M$ be an indecomposable complex of $\ell$-permutation
$(kGb,kHb')$-bimodules.
Then, $M$ induces a stable equivalence between $kGb$ and $kHb'$ if and only if
given $(D,b_D)$ a maximal $b$-subpair, there is 
a maximal $b'$-subpair $(D',b'_{D'})$, an isomorphism
$\phi:D\iso D'$ inducing an isomorphism $\CB r(D,b_D)\iso
\CB r(D',b'_{D'})$ such that
\begin{itemize}
\item
The indecomposable modules occurring in $M$ have vertex-subpairs of the form
$(\Delta_\phi(Q),b_Q\otimes b'_{\phi(Q)})$ for some
$(Q,b_Q)\in \CB r(D,b_D)$, with $(\phi(Q),b'_{\phi(Q)})=\phi(Q,b_Q)$.
\item
For $1\not=Q\le D$, then $b_Q\cdot\Br_{\Delta_\phi Q}M\cdot b'_{\phi(Q)}$
induces
a Rickard equivalence between $kC_G(Q)b_Q$ and $kC_H(Q)b'_{\phi(Q)}$, where
$(Q,b_Q)\in\CB r(D,b_D)$ and $(\phi(Q),b'_{\phi(Q)})=\phi(Q,b_Q)$.
\end{itemize}
\end{thm}

\begin{rem}
In \cite{Ri}, Rickard introduced a notion of splendid equivalences
for principal blocks (complexes of
$\ell$-permutation modules with diagonal vertices), later generalized by
Harris \cite{Ha} and Linckelmann \cite{Li5}. Such equivalences were shown
to induce equivalences for blocks of centralizers. In these approaches, an
isomorphism between the defect groups of the two blocks involved was fixed
a priori and vertex-subpairs were assumed to be ``diagonal'' with respect
to the isomorphism.
Theorem \ref{splendid} shows it is actually easier and more
natural to work with no a priori identification, and the property on
vertex-subpairs is actually automatically satisfied.

The second part of the Theorem (local Rickard equivalences $\Rightarrow$ stable
equivalence) generalizes results of Alperin and Brou\'e and is related
to work of Bouc and Linckelmann.

Finally, a more general theory (terms need not be $\ell$-permutation
modules) has been constructed by Puig (``basic equivalences'') \cite{Pubook}.
\end{rem}

Rickard proposed the following strengthening of Conjecture \ref{ADC}:
\begin{conj}
\label{ADCstrong}
If $D$ is abelian, there is a complex of $\ell$-permutation modules inducing
a Rickard equivalence between
$\CO Gb$ and $\CO N_G(D)c$.
\end{conj}

To the best of my knowledge, in all cases where Conjecture \ref{ADC} is known
to hold, then, Conjecture \ref{ADCstrong} is also known to hold.

Conjecture \ref{ADCstrong} is known to hold when $D$ is cyclic
\cite{Ri1,Li1,papp}. In that
case, one can construct a complex with length $2$, but the longer complex
originally constructed by Rickard might be more natural.
The conjecture holds also when $D\simeq (\BZ/2)^2$ \cite{Ri3,Li2,papp}.
In both cases, the representation type is tame.
Note that there is
no other $\ell$-group $P$ for which Conjecture \ref{ADCstrong} is known to
hold for all $D\simeq P$.

Conjecture \ref{ADCstrong}
holds when $G$ is $\ell$-solvable \cite{Da1,Pu1,HaLi},
when $G$ is a symmetric group or a
general linear group (cf \S \ref{secsl2}; the describing characteristic
case $G=SL_2(\ell^n)$ is solved in \cite{Ok}) and when $G$ is a finite group
of Lie type and $\ell|(q-1)$ (cf \S \ref{secdlADC}).
There are many additional special groups for which the conjecture is known to
hold, cf 
{\rm http://www.maths.bris.ac.uk/$\sim$majcr/adgc/adgc.html}.

\subsubsection{Gluing}
\label{secgluing}
Theorem \ref{splendid} suggests an inductive approach to
Conjecture \ref{ADCstrong}: one should solve the conjecture for local subgroups
(say, $C_G(Q)$, $1\not=Q\le D$) and glue the corresponding Rickard
complexes. This would give rise to a complex inducing a stable equivalence,
leaving us with the core problem of lifting a stable equivalence to a
Rickard equivalence. Unfortunately, complexes are not rigid enough to allow
gluing. This problem can be solved by using complexes endowed with some
extra structure \cite{charl,local}. The idea is to use complexes
that have the properties of chain complexes of simplicial complexes: the
key point is the existence of compatible splittings of the Brauer maps
$M^Q\to M(Q)$. One can build an exact category of $\ell$-permutation modules
with compatible splittings of the Brauer maps. The subcategory of
projective objects turns out to have a very simple description in
terms of sets, and we use only this category. For simplicity, we restrict here
to the case of principal blocks.

\smallskip
Let $G$ be a finite group, $\ell$ a prime number, $k$ an algebraically closed
field of characteristic $\ell$, $b$ the principal block idempotent of $kG$,
$D$ a Sylow $\ell$-subgroup of $G$ and
$c$ the principal block idempotent of $H=N_G(D)$. We assume
$D$ is abelian.
We denote by $Z_\ell(G)$ the 
Sylow $\ell$-subgroup of $Z(G)$ and put $Z=\Delta Z_\ell(G)$.

Let $G'$ be a finite group containing $G$ as a normal
subgroup, let $H'=N_{G'}(D)$ and $F=G'/G\iso H'/H$. We assume $F$ is an
$\ell'$-group,  we put $N=\{(g,h)\in G'\times H^{\prime\opp} |
(gG,hH^\opp)\in \Delta F\}$ and
$\bar{N}=N/Z$.

Let $\CE$ be the category of $\bar{N}$-sets whose point stabilizers
are contained in $\Delta D/Z$. Let $\tilde{\CE}$ be the Karoubian envelop
of the linearization of $\CE$ (objects are pairs $(\Omega,e)$ where
$\Omega$ is a $\bar{N}$-set and $e$ an idempotent of the monoid
algebra of $\End_{\bar{N}}(\Omega)$).
We have a faithful functor $\tilde{\CE}\to k\bar{N}\mlperm,\ 
(\Omega,e)\mapsto k(\Omega,e):=k\Omega e$.

\smallskip
We are now ready to state a further strengthening of Conjecture \ref{ADC}.
For the inductive approach, it is important to take into account central
$\ell$-subgroups and $\ell'$-automorphism groups.

\begin{conj}
\label{ADCstrong2}
There is a complex $C$ of objects of $\tilde{\CE}$ such that
$\Res_{G\times H^\opp}^{\bar{N}}k(C)$ induces a Rickard equivalence between
$kGb$ and $kHc$.
\end{conj}

We can also state a version of Question \ref{liftability}, for the pair
$(G',G)$:
\begin{question}
\label{conjlift}
Let $C$ be a complex of objects of $\tilde{\CE}$ such that
$\Res_{G\times H^\opp}^{\bar{N}}k(C)$ induces a stable equivalence between
$kGb$ and $kHc$.
Is there a bounded complex $R$ of finitely generated
projective $\bar{N}$-modules and a morphism $f:R\to k(C)$ such that
$\Res_{G\times H^\opp}^{\bar{N}}\cone(f)$  induces a Rickard equivalence
between $kGb$ and $kHc$?
\end{question}

The following Theorem reduces (a suitable version of)
the abelian defect conjecture to (a suitable version of)
the problem of lifting stable equivalences to Rickard equivalences.

\begin{thm}
Assume Question \ref{conjlift} has a positive answer for
$(N_{G'}(Q),C_G(Q))$ for all non trivial subgroups $Q$ of $D$. Then,
Conjecture \ref{ADCstrong2} holds.
\end{thm}

The proof goes by building inductively (on the index of $Q$ in $D$) a system
of complexes for $N_{G'}(Q)$ and gluing them together. The key point
is that, given a finite group $\Gamma$, the category of $\Gamma$-sets whose
point stabilizers are non-trivial $p$-subgroups is locally determined. This
allows us to manipulate objects of $\tilde{\CE}$ as ``sheaves''.

\subsection{Chevalley groups}
\label{secChev}
We explain Brou\'e's idea that complexes of cohomology of certain varieties
should give rise to derived equivalences, for finite groups of Lie type.

\subsubsection{Deligne-Lusztig varieties}
\label{secDL}
Let $\BG$ be a connected reductive algebraic group defined over a finite
field and let $F$ be an endomorphism of $\BG$, a power $F^d$ of which is a
Frobenius endomorphism defining a structure over a finite field $\BF_{q^d}$
for some $q\in \BR_{>0}$. Let $G=\BG^F$ be the associated finite group.

Let $\ell$ be a prime number with $\ell {\not|}\ q$,
$K$ a finite extension of $\BQ_l$, and $\CO$ its ring of integers. We assume
$K$ is big enough.

\smallskip
Let $\BL$ be an $F$-stable Levi subgroup of $\BG$,
$\BP$ be a parabolic subgroup with Levi complement
$\BL$, and let $\BU$ be the unipotent radical of $\BP$. We define the
Deligne-Lusztig variety
$$Y_\BU=\{g\BU\in \BG/\BU \mid g^{-1}F(g)\in \BU\cdot F(\BU)\},$$
a smooth affine variety with a left action of $\BG^F$ and a right
action of $\BL^F$ by multiplication. The corresponding complex
of cohomology $R\Gamma_c(Y_\BU,\CO)$ induces the
Deligne-Lusztig induction functor $R_{\BL\subset\BP}^\BG:
D^b(\CO \BL^F)\to D^b(\CO \BG^F)$.

The effect of these functors on characters (\ie, $K_0$'s after extension to
$K$) is a central tool for Deligne-Lusztig and Lusztig's construction of
irreducible characters of $G$.
It is important to also consider the finer invariant $\tR \Gamma_c(Y_\BU,\CO)$,
an object of $K^b(\CO(\BG^F\times (\BL^F)^\opp)\mlperm)$
which is quasi-isomorphic to $R\Gamma_c(Y_\BU,\CO)$
\cite{Rietale,affine}.

We put $X_\BU=Y_\BU/\BL^F$ and denote by $\pi:Y_\BU\to X_\BU$ the quotient
map.

\begin{rem}
One could use ordinary cohomology instead of the compact support version.
One can conjecture that the two versions are interchanged by
Alvis-Curtis duality:
$(R\Gamma_c(Y_\BU,\CO)\otimes_{\CO \BL^F}^\BL -)\circ D_\BL$ and
$D_\BG \circ (R\Gamma(Y_\BU,\CO))\otimes_{\CO \BL^F}^\BL -)$ should
differ by a shift.
This is known in the Harish-Chandra case, \ie, when $\BP$ is $F$-stable 
\cite{CabRi}.
\end{rem}

\smallskip
Let $\BT_0\subset\BB_0$ be a pair consisting of an $F$-stable maximal torus
and an $F$-stable Borel subgroup of $\BG$. Let $\BU_0$ be the unipotent
radical of $\BB_0$ and let $W=N_\BG(\BT_0)/\BT_0$.

Let $B^+$ (resp. $B$) be the braid monoid (resp. group) of $W$.
The canonical map
$B^+\to W$ has a unique section $w\mapsto \Bw$ that preserves lengths
(it is not a group morphism!).
We fix an $F$-equivariant
morphism $\tau:B\to N_\BG(\BT_0)$ that lifts the canonical map
$N_\BG(\BT_0)\to W$ \cite{Ti}.
Given $w\in W$, we put $\dot{w}=\tau(\Bw)$.
Let $w_0$ be the longest element of $W$ and let
$\Bpi=\Bw_0^2$, a central element of $B$.

\smallskip
Assume $\BL$ above is a torus. We give a different model for $Y_\BU$.
Let $w\in W$ and $h\in\BG$ such that $h^{-1}F(h)=\dot{w}$ and
$\BU=h\BU_0 h^{-1}$. Let
$$Y(w)=\{g\BU_0\in \BG/\BU_0 \mid g^{-1}F(g)\in \BU_0 \dot{w}\BU_0\},$$
a variety with a left action of $G$ and a right action of
$\BT_0^{wF}$ by multiplication.
We have $\BL=h\BT_0 h^{-1}$ and conjugation by $h$ induces an isomorphism
$\BL^F \iso \BT_0^{wF}$.
Right multiplication by $h$ induces an isomorphism
$Y_\BU\iso Y(w)$ compatible with the actions of $G$ and $\BL^F$.
We have $\dim Y(w)=l(w)$. We write
$Y_F(w)$ when the choice of $F$ is important.

\smallskip
Given $w_1,\ldots,w_r\in W$, we put
\begin{multline*}
Y(w_1,\ldots,w_r)=\{(g_1\BU_0,\ldots,g_r\BU_0)\in (\BG/\BU_0)^r
\mid\\
g_1^{-1}g_2\in \BU_0 \dot{w}_1 \BU_0,\ldots,
g_{r-1}^{-1}g_r\in \BU_0 \dot{w}_{r-1} \BU_0 \textrm{ and }
g_r^{-1}F(g_1)\in \BU_0 \dot{w}_r \BU_0\}.
\end{multline*}

Up to a transitive system of canonical isomorphisms,
$Y(w_1,\ldots,w_r)$ depends only on the product $b=\Bw_1\cdots\Bw_r\in
B^+$ and we denote that variety by $Y(b)$ \cite{De,BrMi2}.

\subsubsection{Jordan decomposition}
\label{secJordan}
As a first step in his
classification of (complex) irreducible characters of finite groups of
Lie type, Lusztig established a
{\em Jordan decomposition} of characters.

Let $(\BG^*,F^*)$ be Langlands dual to $(\BG,F)$. Then, Lusztig defined
a partition of the set $\Irr(G)$ of irreducible characters of $G$:
$$\Irr(G)=\coprod_{(s)} \Irr(G,(s))$$
where $(s)$ runs over conjugacy classes of semi-simple elements of
$(\BG^*)^{F^*}$.
The elements in $\Irr(G,1)$ are the {\em unipotent characters}.

Furthermore, Lusztig constructed a bijection
\begin{equation}
\label{decJordan}
\Irr(G,(s))\iso \Irr((C_{\BG^*}(s)^*)^F,1)
\end{equation}
(assuming $C_{\BG^*}(s)$ is connected).
So, an irreducible character corresponds to a pair consisting of a semi-simple
element in the dual and a unipotent character of the dual of the
centralizer of that semi-simple element.

\medskip
Brou\'e and Michel \cite{BrMi}
showed that the union of series corresponding to classes
with a fixed $\ell'$-part is a union of blocks:
let $t$ be a $\ell'$-element of $(\BG^*)^{F^*}$ and let
$$\Irr(G,(t))_{\ell}=\coprod_{(s)}\Irr(G,(s))$$
where $(s)$ runs over conjugacy classes of semi-simple elements
of $(\BG^*)^{F^*}$ whose $\ell'$-part is conjugate to
$t$. Then, $\Irr(G,(t))_{\ell}$ is a union of $\ell$-blocks and
we denote by $B(G^F,(t))$ the corresponding factor of $\CO G^F$.

Brou\'e \cite{BrIHES} conjectured that the decomposition (\ref{decJordan})
arises from a Morita equivalence (cf also \cite{Hi}). More, precisely,
we have the following Theorem \cite[Theorem B']{BoRou1} obtained in joint
work with
C.~Bonnaf\'e (cf also \cite{CabEn} for a detailed exposition).
This was conjectured by Brou\'e who gave a proof when
$t$ is regular \cite{BrIHES}.

\begin{thm}[Jordan decomposition of blocks]
\label{jordan}
Assume $C_{\BG^*}(t)$ is contained in an $F^*$-stable
Levi subgroup $\BL^*$ of $\BG^*$
with dual $\BL\le\BG$.
Let $\BP$ be a parabolic subgroup of $\BG$ with Levi
complement $\BL$ and unipotent radical $\BU$. Let $d=\dim X_\BU$ and let
$\CF_t=\pi_*\CO\otimes_{\CO\BL^F}B(\BL^F,(t))$.

Then,
$H^i_c(X_\BU,\CF_t)=0$ for $i\not=d$ and
$H^d_c(X_\BU,\CF_t)$
induces a Morita equivalence between  
$B(G,(t))$ and $B(\BL^F,(t))$.
\end{thm}

The Theorem reduces the study of blocks of finite groups of Lie type to
the case of those associated to a quasi-isolated element $t$. When
$\BL^*=C_{\BG^*}(t)$ is a Levi subgroup of $\BG^*$, then
$B(\BL^F,(t))$ is isomorphic to $B(\BL^F,1)$.

As shown by Brou\'e, the key point is the statement about the vanishing of
cohomology.
When $\BL$ is a torus, this
is \cite[Theorem 9.8]{DeLu}. For the general case,
two
difficulties arise: there are no known good smooth compactifications of
the varieties $X_\BU$ and the locally constant sheaf $\CF_t$
has wild ramification.
We solve these issues as follows. Let $\bar{X}$ be the closure of $X_\BU$ in
$\BG/\BP$. We construct new varieties of Deligne-Lusztig type and commutative
diagrams
$$\xymatrix{
X_i\ar@{^{(}->}[r]^{j_i}\ar[d]_{f'_i} & Y_i \ar[d]^{f_i} \\
X_\BU\ar@{^{(}->}[r]_{j} & \bar{X}
}$$
where $Y_i$ is smooth, $Y_i-X_i$ is a divisor with normal crossings, and
$f_i$ is proper. We also construct tamely ramified sheaves $\CF_i$ on $X_i$
with the following properties:
\begin{itemize}
\item $\CF_t$ is in the thick subcategory of the derived category of
constructible sheaves on $\BX_\BU$ generated by the $Rf'_{i*}\CF_i$
\item $(Rj_{i*}\CF_i)_{|f_i^{-1}(\bar{X}-X_\BU)}=0$.
\end{itemize}

The first property follows from the following generation result of the derived
category of a finite group of Lie type \cite[Theorem A]{BoRou1}:

\begin{thm}
The category of perfect complexes for $B(G,(t))$ is generated, as a thick
subcategory, by the $R_{\BT\subset\BB}^\BG B(\BT^F,(t))$, where
$\BT$ runs over the $F$-stable maximal tori of $\BG$ such that $t\in\BT^*$
and $\BB$ runs over the Borel subgroups of $\BG$ containing $\BT$.
\end{thm}

\begin{rem}
Note that the corresponding result for derived categories is true, under
additional assumptions on $\BG$ \cite{BoRou3}:
this is related to Quillen's Theorem,
we need every elementary abelian $\ell$-subgroup of $G$ to be contained
in an $F$-stable torus of $\BG$.
\end{rem}

\begin{rem}
Note that
the Morita equivalence of Theorem \ref{jordan} is not splendid in general.
This issue is analyzed in \cite{BoRou3}.
\end{rem}

\begin{example}
Let 
$\BG=\GL_n(\bar{\BF}_q)$ and $F:(x_{ij})_{1\le i,j\le n}\mapsto 
(x_{ij}^q)_{i,j}$. We have $G=\GL_n(\BF_q)$, $\BG=\BG^*$ and $F^*=F$.
Centralizers of semi-simple elements are Levi subgroups, so Theorem
\ref{jordan} gives a Morita equivalence between any block of a
general linear group over $\CO$ and a unipotent block.

\end{example}

\subsubsection{Abelian defect conjecture}
\label{secdlADC}
Let $b$ be a block idempotent
of $\CO G$. Let $(D,b_D)$ be a maximal $b$-subpair, let $H=N_G(D,b_D)$
and let $\BL=C_\BG(D)$.
We assume $D$ is abelian and $\BL$ is a Levi subgroup of $\BG$ (these are
satisfied if $\ell{\not|}\ |W|$).

Brou\'e conjectured that the sought-for complex in Conjecture \ref{ADCstrong}
should arise from Deligne-Lusztig varieties (\cite[p.81]{BrLuminy},
\cite[\S 1]{BrMa}, \cite[\S VI]{Brmit}):

\begin{conj}[Brou\'e]
\label{ADCdlvar}
There is a parabolic subgroup $\BP$ of $\BG$ with Levi complement $\BL$
and unipotent radical $\BU$,
and a complex $C$ inducing a Rickard equivalence between $\CO Gb$ and
$\CO H b_D$ such that $\Res_{G\times (\BL^F)^\opp}C$ is isomorphic
to $\tR\Gamma_c(Y_\BU,\CO)b_D$.
\end{conj}

This conjecture \ref{ADCdlvar} is known to hold \cite{PuLuminy}
when there is a choice of an $F$-stable
parabolic subgroup $\BP$ (case $\ell| (q-1)$).
Then, $Y_\BU$ is $0$-dimensional and
the Deligne-Lusztig induction is the Harish-Chandra induction. The key steps in
the proof are:
\begin{itemize}
\item Produce an action of the reflection group
$H/\BL^F$ from a natural action of the associated Hecke algebra. One needs
to show that certain obstructions vanish.
\item Identify a $2$-cocycle of $H/\BL^F$ with values in $\CO^\times$.
\item Compute the dimension of the $KG$-endomorphism ring.
\end{itemize}

\subsubsection{Regular elements}
\label{secreg}
As a first step, one should make
Conjecture  \ref{ADCdlvar} more precise by specifying $\BP$ and by defining
the extension of the action of $C_G(D)$ to an action of $H$
on $\tilde{R}\Gamma_c(Y_\BU,\CO)b_D$. These issues are partly solved and I will
explain
the best understood case where $\BL=\BT$ is a torus and $\CO Gb$ is the
principal block (cf \cite{BrMi2}).
Assume as well $\ell{\not|}\ (q-1)$. To simplify, assume
further that $F$ acts trivially on $W$ (``split'' case).

Note that $\BT$ defines a conjugacy class $\CC$ of $W$ and
the choice of $\BP$ amounts to the choice of $w\in \CC$
(defined from $\BP$ as in \S \ref{secDL}).
Since $\BT=C_\BG(D)$, it follows that elements in $\CC$ are
Springer-regular.
There is $w_d\in \CC$ such that $(\Bw_d)^d=\Bpi$, where
$d>1$ is the order of $w_d$  (a ``good'' regular element).

\smallskip
Given $w\in W$, we have a purely inseparable morphism
\begin{align*}
Y(w,w^{-1}w_0,w_0ww_0,w_0w^{-1})&\to Y(w^{-1}w_0,w_0ww_0,w_0w^{-1},w)\\
(x_1,x_2,x_3,x_4)&\mapsto (x_2,x_3,x_4,F(x_1)).
\end{align*}
Via the canonical isomorphisms, this induces an endomorphism of $Y(\Bpi)$.
This extends to an action of $B^+$ on $Y(\Bpi)$.

There is an embedding of $Y_F(w_d)$ as a closed subvariety of
$Y_{F^d}(w_d,\ldots,w_d)$ ($d$ terms) given by
$$x\mapsto (x, F(x),\ldots,F^{d-1}(x)).$$
The action of $C_{B^+}(\Bw_d)$ on $Y_{F^d}(\Bpi)$ restricts to an action
on $Y_F(w_d)$. It induces an action of $C_B(\Bw^d)$ on
$\tilde{R}\Gamma_c(Y(\Bw_d),\CO)$.

The group $H/C_G(D)\simeq C_W(w_d)$ is a complex reflection group and we
denote by $B_d$ its braid group. There
is a morphism $B_d\to C_B(\Bw_d)$,
uniquely defined up to conjugation by an element
of the pure braid group of $C_W(w_d)$ (it is expected to be an isomorphism,
and known to be such in a number of cases \cite{BesDiMi}).

Now, the conjecture is that, up to homotopy, the action of
$\CO(\BT_0^{w_d F}\rtimes B_d)$ on $\tR\Gamma_c(Y(w_d),\CO)b_D$ induces an
action of the quotient algebra $\CO Hc$ and the resulting object is a splendid
Rickard complex:

\begin{conj}
\label{ADCdlvarstrong}
There is a complex $C\in K^b((\CO Gb)\otimes(\CO Hb_D)^\opp\mlperm)$, unique up
to isomorphism, with the following properties:
\begin{itemize}
\item There is a surjective morphism
$f:\CO \BT_0^{w_d F}\rtimes B_d\to \CO Hb_D$ extending the inclusion
$\BT_0^{w_d F}\subset H$ such that
\begin{itemize}
\item $f^*C$ and $\tR\Gamma_c(Y(w_d),\CO)b_D$
are isomorphic in 
$D^b((\CO \BT_0^{w_d F}\rtimes B_d)\otimes(\CO H)^\opp)$
\item the map $kB_d\to kC_W(w_d)$ deduced from $f$ by applying
$k\otimes_{\CO \BT_0^{w_d F}}-$ is the canonical map.
\end{itemize}
\item $C$ is isomorphic to
$\tR\Gamma_c(Y(w_d),\CO)b_D$
in $K^b((\CO G)\otimes(\CO C_G(D))^\opp\mlperm)$.
\end{itemize}
Furthermore, such a complex $C$ induces
a Rickard equivalence between $\CO Gb$ and $\CO Hb_D$.
\end{conj}

The most crucial and difficult part in that conjecture is to show that we
have no non-zero shifted endomorphisms of the 
complex (``disjunction property''),
either for the action of $G$ or for that $H$.

\smallskip
Conjecture \ref{ADCdlvarstrong} is known to hold when $l(w_d)=1$
\cite{affine} and for $\GL_n$ and $d=n$ \cite{BoRou2}. In the first case,
we use good properties of cohomology of curves and prove disjunction for
the action of $G$. In the second case, we study the variety 
$D(\BU_0)^F\backslash Y(w_d)$ and prove disjunction for the action of $H$.
This works only for $\GL_n$, for we 
rely on the fact that induced Gelfand-Graev representations generate
the category of projective modules.

\begin{rem}
The version ``over $K$'' of Conjecture \ref{ADCdlvarstrong} is open, even
after restricting to unipotent representations 
(=applying the functor $K\otimes_{KT^{w_dF}}-$). The action of
$KB_d$ on $H_c^*(Y(w_d),K)$ should factor through an action of the Hecke
algebra of $C_W(w_d)$, for certain parameters. This is known in some
cases, in particular when $d=2$ (work of Lusztig \cite{LuCBMS}
and joint work with Digne and Michel \cite{DiMiRou}) and for
other cases \cite{DiMi}. The disjunction property is
known for $w_d$ a Coxeter element \cite{Lucox}, for $\GL_n$ and $d=n-1$
\cite{DiMi} and in most rank $2$ groups \cite{DiMiRou}.
\end{rem}

\begin{rem}
When $\ell|(q-1)$ (case $d=1$), one can formulate a version of Conjecture 
\ref{ADCdlvarstrong} using the variety $Y(\Bpi)$ \cite[Conjectures 2.15]{BrMi2}.
\end{rem}

\subsection{Local representation theory as non-commutative birational
geometry}
\label{secgeo}
It is expected that birational Calabi-Yau varieties should have equivalent
derived categories (cf \cite{Bri2}). We view Question \ref{liftability} as a
non-commutative version: one can expect that ``sufficiently nice''
Calabi-Yau triangulated categories are determined by (not too small) quotients.
We explain here how this analogy can be made precise, in the setting of
McKay's correspondence, via Koszul duality.

\subsubsection{$2$-elementary abelian defect groups}
Let $P$ be an elementary abelian $2$-group. Let $k$ be a field of
characteristic $2$ and $V=P\otimes_{\BF_2}k$. Let $E$ be a group of odd
order of automorphisms of $P$. The algebras $kP\rtimes E$ and
$\Lambda(V)\rtimes E$ are isomorphic.

Koszul duality (cf eg \cite{KeDG}) gives an equivalence
$$D^b((\Lambda(V)\rtimes E)\mModgr)\iso
D^b_{E\times \BG_m}(V).$$

\subsubsection{McKay's correspondence}
Let $V$ be a finite-dimensional vector space over $k$ and $E$ a finite
subgroup of $\GL(V)$ of order invertible in $k$. Recall the following
conjecture (independence of the crepant resolution):

\begin{conj}[McKay's correspondence]
\label{conjMcKay}
If $X\to V/E$ is a crepant resolution, then $D^b(X)\simeq D^b_E(V)$.
\end{conj}

The conjecture is known to hold when $\dim V=3$ \cite{BriKiRei,Bri1} (in
dimension $3$, the Hilbert scheme of $E$-clusters on $V$ is a crepant
resolution).
It is also known when $V$ is a symplectic vector space and $E$ respects the
symplectic structure \cite{BezKa}. See \cite[\S 2.2]{Bri2} for more details.

\smallskip
Examples in dimension $>3$ where $E-\Hilb V$ is smooth are rare.
An infinite family of examples is provided by the following Theorem of
Sebestean \cite{Se}:

\begin{thm}
\label{magda}
Let $n\ge 2$, let $k$ be a field containing a primitive
$(2^n-1)$-th root of unity $\zeta$ and let
$E$ be the subgroup of $\SL_n(k)$ generated by the diagonal matrix with
entries
$(\zeta,\zeta^2,\ldots,\zeta^{2^{n-1}})$. Assume $2^n-1$ is invertible in $k$.

Then, $E-\Hilb(\BA_k^n)$ is a smooth crepant resolution of
$\BA_k^n/E$ and there is an equivalence
$D^b_E(\BA_k^n)\iso D^b(E-\Hilb(\BA_k^n))$.
\end{thm}
The diagonal action of $\BG_m$ on $\BA_k^n$ induces an action on
$E-\!\Hilb(\BA_k^n)$ and the equivalence is equivariant for these actions.

\smallskip
Let $G=\SL_2(2^n)$, let $P$ be the subgroup of strict upper triangular matrices
(a Sylow $2$-subgroup), and let $E$ be the subgroup of diagonal matrices.
The action of $E$ on $P\otimes_{\BF_2}\bar{\BF}_2$ coincides with the one in
Theorem \ref{magda}.
Combining the solution of Conjecture
\ref{ADC} for $G$ (Okuyama, \cite{Ok}) and \S \ref{secgradabelian},
 the Koszul duality equivalence, and Theorem \ref{magda},
 we deduce a geometric realization of modular representations of
$\SL_2(2^n)$ in natural characteristic:

\begin{cor}
There is a grading on the principal $2$-block $A$
of $\bar{\BF}_2G$ and an equivalence
$D^b(A\mModgr)\iso D^b_{\BG_m}(E-\Hilb \BA^n_k)$.
\end{cor}

\begin{rem}
It should be interesting to study homotopy categories of sheaves on
singular varieties and their relation to derived categories of crepant
resolutions.
\end{rem}

\subsection{Perverse Morita equivalences}
\label{secperv}
In this part, we describe joint work with J.~Chuang \cite{CY}.

\subsubsection{Definitions}
Let $\CA,\CA'$ be two abelian categories. We assume every
object has a finite composition series.
Let $\CS$ (resp. $\CS'$) be the set of isomorphism classes of simple objects
of $\CA$ (resp. $\CA'$).

\begin{defi}
An equivalence $F:D^b(\CA)\iso D^b(\CA')$ is perverse if there is
\begin{itemize}
\item a filtration
$\emptyset=\CS_0\subset \CS_1\subset\cdots\subset \CS_r=\CS$
\item a filtration
$\emptyset=\CS'_0\subset \CS'_1\subset\cdots\subset \CS'_r=\CS'$
\item and a function $p:\{1,\ldots,r\}\to\BZ$
\end{itemize}
 such that
\begin{itemize}
\item $F$ restricts to equivalences $D^b_{\CA_i}(\CA)\iso D^b_{\CA'_i}(\CA')$
\item $F[-p(i)]$ induces equivalences $\CA_i/\CA_{i-1}\iso
\CA'_i/\CA'_{i-1}$.
\end{itemize}
where $\CA_i$ (resp. $\CA'_i$) is the Serre subcategory of $\CA$
(resp. $\CA'$) generated by $\CS_i$ (resp. $\CS'_i$).
\end{defi}

An important point is that $\CA'$ is determined, up to equivalence,
by $\CA$, $\CS_\bullet$ and $p$.

\subsubsection{Symmetric algebras}
Let $A$ be a symmetric finite dimensional algebra and $\CA=A\mMod$.
Then, given $\CS_\bullet$ and $p$, one can construct
a tilting complex $T$ with $A'=\End_{D^b(A)}(T)$, $\CA'=A'\mMod$
 and $F=R\Hom_A^\bullet(T,-)$
(this cannot be achieved in general for a non-symmetric algebra $A$).

\begin{rem}
One might ask whether all derived equivalences between finite dimensional
symmetric algebras are compositions of perverse equivalences, or at least,
if two derived equivalent
symmetric algebras can be related by a sequence of perverse equivalences.
The known examples in block theory seem to all be composition of perverse
equivalences.
\end{rem}

\begin{rem}
One can expect the equivalences predicted in Conjecture \ref{ADCdlvar}
will be perverse. The filtration should be provided by Lusztig's $a$-function.
\end{rem}

Iterating the construction of tilting complexes associated to
$\CS_\bullet$ and $p$, one gets
an action of $\operatorname{Free}(\CP(\CS))\rtimes\GS(\CS)$ on
the set of $(T_s)_{s\in\CS}$, where $T_s$ is an indecomposable perfect
complex of $A$ and $\bigoplus_s T_s$ is tilting. We expect this will
relates to Bridgeland's space of stability conditions
\cite[\S 4]{Bri2} (such spaces seem to have a nicer behaviour under the
Calabi-Yau assumption).

\begin{rem}
The considerations above are interesting for Calabi-Yau algebras of positive
dimension, though, in that case, tilting is not always possible. When
$r=2$ and $|\CS_2-\CS_1|=1$, tilting has been known in string theory
as Seiberg duality.
\end{rem}

\section{Invariants}
\label{secinvar}
Invariants of triangulated categories and dg-categories are
discussed in \cite[\S 6]{KeICM}. We discuss here some more elementary
invariants,
used to study finite dimensional algebras.

\subsection{Automorphisms of triangulated categories}
\label{secauto}
\subsubsection{Rings}

Let $k$ be a commutative ring and $A$ be a $k$-algebra. We denote by
$\Pic(A)$ the group of isomorphism classes of
invertible $(A,A)$-bimodules and by
$\DPic(A)$ the group of isomorphism classes of invertible objects of the
derived category of $(A,A)$-bimodules: this is the part of the automorphism
group of $D(A\mMOD)$ that comes from standard equivalences. By Rickard's
Theorem, $\DPic(A)$ is invariant under derived equivalences.

The following Proposition
has been observed by many people (Rickard, Roggenkamp-Zimmermann,
\cite[Proposition 3.3]{RouZi}, \cite[Proposition 3.4]{Ye1},...).

\begin{prop}
If $A$ is local, then $\DPic(A)=\Pic(A)\times
\langle A[1]\rangle$.
\end{prop}

Given $R$ a flat commutative $Z$-algebra, there is a canonical
morphism $\DPic(A)\to \DPic(A\otimes_Z R)$ (joint work with
A.~Zimmermann \cite[\S 2.4]{RouZi}). If
$R$ is faithfully flat over $Z$, the kernel of that map is contained in
$\Pic(A)$. This is the key point for the following (cf
\cite[Proposition 3.5]{Ye1} and \cite[Proposition 3.3]{RouZi}):

\begin{thm}
Assume $A$ is commutative and indecomposable. Then, $\DPic(A)=\Pic(A)\times
\langle A[1]\rangle$.
\end{thm}

\subsubsection{Invariance of automorphisms}
Let $A$ be a finite dimensional algebra over an algebraically closed field
$k$. We denote by $\Aut(A)$ the group of
automorphisms of $A$. This is an algebraic group and we denote by
$\Inn(A)$ its closed subgroup of inner automorphisms. We put
$\Out(A)=\Aut(A)/\Inn(A)$. We have a morphism of groups
$\Aut(A)\to\Pic(A),\ \alpha\mapsto [A_\alpha]$, where
$A_\alpha=A$ as a left $A$-module and the right action of $a\in A$ is
given by right multiplication by $\alpha(a)$. It induces an injective
morphism $\Out(A)\to \Pic(A)$.

\smallskip
The following result \cite{grad}
gives a functorial interpretation of $\Out$, to be compared with the
functorial interpretation of $\Pic(X)$ for a smooth projective variety $X$.

\begin{thm}
\label{functorOut}
The functor from the category of affine varieties over $k$ to groups that
sends $X$ to the set of isomorphism classes of
$(A\otimes A^\opp\otimes\CO_X)$-modules
that are locally free of rank $1$ as $(A\otimes\CO_X)$ and as
$(A^\opp\otimes\CO_X)$-modules is represented by $\Out(A)$.
\end{thm}

The following Theorem \cite{grad}
shows the invariance of $\Out^0$, the identity component
of $\Out$, under certain
equivalences. In the case of Morita equivalences, it goes back to Brauer, and
for derived equivalences, it
has been obtained independently by Huisgen-Zimmermann and Saor\'{\i}n
\cite{HuiSa}. In these cases, it follows easily from Theorem \ref{functorOut}
while, for stable equivalences, some work is needed to get rid globally of
projective direct summands.

\begin{thm}
\label{Out}
Let $B$ be a finite dimensional $k$-algebra and
let $C$ be a bounded complex of finitely generated
$(A,B)$-bimodules inducing a derived equivalence
or a stable equivalence (in which case we assume $A$ and $B$ are
self-injective). Then, there is a unique
isomorphism of algebraic groups
$\sigma:\Out^0(A)\iso \Out^0(B)$ such that
$A_\alpha\otimes_A C\simeq C\otimes_B B_{\sigma(\alpha)}$ for all
$\alpha\in\Out^0(A)$.
\end{thm}

Yekutieli \cite{Ye2} deduces that $\DPic(A)$ has a structure of a locally 
algebraic group, with connected component $\Out^0(A)$.

\subsubsection{Coherent sheaves}

The following result \cite{grad} is a variant of Theorem \ref{Out}.
\begin{thm}
Let $X$ and $Y$ be two smooth projective schemes over an algebraically closed
field $k$. An equivalence $D^b(X)\iso D^b(Y)$ induces an
isomorphism $\Pic^0(X)\rtimes\Aut^0(X)\iso \Pic^0(Y)\rtimes \Aut^0(Y)$.
\end{thm}

This implies in particular that if $A$ and $B$ are derived equivalent
abelian varieties, then there is a symplectic isomorphism
$\hat{A}\times A\iso \hat{B}\times B$ (and the converse holds as well
\cite{Or,Po}).

\subsubsection{Automorphisms of stable categories and endo-trivial modules}

Let $A$ be a finite dimensional self-injective algebra over an algebraically
closed field $k$. We denote by $\StPic(A)$ the group of isomorphism classes
of invertible objects of $(A\otimes A^\opp)\mstab$.

\smallskip
Let $P$ be an $\ell$-group and $k$ a field of characteristic $\ell$. A finitely
generated $kP$-module $L$ is an {\em endo-trivial} module if
$L\otimes_k L^*\simeq k$ in $kP\mstab$ or equivalently, if
$\End_{kP\mstab}(L)=k$ \cite{Car}.
Note that the classification of endo-trivial modules has been recently
completed \cite{CarTh} (the case where $P$ is abelian goes back to \cite{Da2}).

Let $\CT(kP)$ be the group of isomorphism classes of indecomposable
endo-trivial modules.
We have an injective morphism of groups
$$\CT(kP)\to \StPic(kP),\ [L]\mapsto
[\Ind_{\Delta P}^{P\times P^\opp}L].$$
This extends to an isomorphism $\CT(kP)\times \Out(kP)\iso \StPic(kP)$
(\cite[\S 3]{Listab} and \cite[\S 2]{CarRou}). 

Let $Q$ be an $\ell$-group.
A stable equivalence of Morita type $kP\mstab\iso kQ\mstab$ induces
an isomorphism $\CT(kP)\iso \CT(kQ)$.
It actually forces the algebras
$kP$ and $kQ$ to be isomorphic
(\cite[\S 3]{Listab}, \cite[Corollary 2.4]{CarRou}).
It is an open question whether this implies that $P$ and $Q$ are isomorphic.

\begin{thm}[{\cite[Theorem 3.2]{CarRou}}]
\label{liftFrob}
Let $P$ be an abelian $\ell$-group and $E$ a cyclic $\ell'$-group acting 
freely on $P$. We put $G=P\rtimes E$. Then,
$\StPic(kG)=\Pic(kG)\cdot \langle \Omega\rangle$.
In particular, the canonical morphism
$\TrPic(kG)\to \StPic(kG)$ is surjective.
\end{thm}

\begin{rem}
Let $A$ be a block over $k$ of a finite group, with defect group isomorphic to
$P$ and $N_G(P)/P$ acting as $E$ on $P$.
From Theorem \ref{liftFrob}, one deduces \cite[Corollary 4.4]{CarRou} via a
construction of Puig \cite{Pugeom}, that
a stable equivalence of Morita
type between $A$ and $kG$ lifts to a Rickard equivalence if and only if
$A$ and $kG$ are Rickard equivalent if and only if they are splendidly
Rickard equivalent. In particular, for blocks with abelian defect group
$D$ such that $N_G(D,b_D)/C_G(D)$ is cyclic, then Conjecture \ref{ADC} implies
Conjecture \ref{ADCstrong}.
\end{rem}

\subsection{Gradings}
\label{secgrad}
In this section, we describe results of \cite{grad}.

\subsubsection{Transfer of gradings}

We assume we are in the situation of Theorem \ref{Out}.
Assume $A$ is graded, \ie, there is a morphism $\BG_m\to\Aut(A)$. The
induced morphism $\BG_m\to\Out^0(A)$ induces a morphism
$\BG_m\to\Out^0(B)$. There exists a lift to a morphism
$\BG_m\to \Aut^0(B)$, and this corresponds to a grading on $B$.
There is a grading on (an object isomorphic to) 
$C$ that makes it into a complex of graded
$(A,B)$-bimodules and it induces an equivalence between the appropriate
graded categories.

Let $A$ be a self-injective indecomposable graded algebra, let $n$ be the
largest integer such that $A_n\not=0$, and let
$C\in\BZ[q,q^{-1}]$ be the graded Cartan matrix of $A$.

If $A$ is non-negatively
graded and the Cartan matrix of $A_0$ has non-zero determinant, then
$\deg \det(C)=nr$, where $r$ is the number of simple $A$-modules. As a
consequence, one gets a positive solution of a ``non-negatively graded''
version of Conjecture \ref{stablesimple}:

\begin{prop}
Let $A$ and $B$ be two indecomposable self-injective non-negatively
graded algebras. Assume $A_0$ has finite global dimension and
there is a graded stable equivalence of Morita type
between $A$ and $B$.
Then, $A$ and $B$ have the same number of
simple modules.
\end{prop}

\begin{rem}
Let $A$ be a non-negatively graded indecomposable self-injective algebra
with $A_0$ of finite global dimension. Let $B$ be a stably equivalent
self-injective algebra. One could hope that there is a compatible grading on
$B$ that is non-negative, but this is not possible in general.
It would be still be very interesting to see if this can be achieved
if the grading on $A$ is ``tight'' in the sense of Cline-Parshall-Scott,
\ie, if $\bigoplus_{j\le i}A_j=(JA)^i$ (cf the gradings in
\S \ref{secgradabelian}).
\end{rem}

\subsubsection{Blocks with abelian defect}
\label{secgradabelian}
Let $P$ be an abelian $\ell$-group and $k$ an
algebraically closed field of characteristic $\ell$.
The algebra $kP$ is (non-canonically) isomorphic to the graded algebra
associated to the radical filtration of $kP$. Fixing such an isomorphism
provides a grading on $kP$. Let $E$ be an $\ell'$-group of automorphisms of
$P$. Then, the isomorphism above can be made $E$-equivariant and we
obtain a structure of graded algebra on $kP\rtimes E$ extending the
grading on $kP$ and with $kE$ in degree $0$. Given a central extension of
$E$ by $k^\times$, this construction applies as well
to the twisted group algebra $k_* P\rtimes E$.

Let $A$ be a block of a finite group over $k$
with defect group $D$. Then, there is $E$ and
a central extension as above such that the corresponding block of $N_G(D)$
is Morita equivalent to $k_* D\rtimes E$ \cite{Ku}.
So, Conjecture \ref{ADC} predicts there are interesting gradings on $A$. In
the inductive approach to Conjecture \ref{ADCstrong2}, there is a stable
equivalence of Morita type between $A$ and $k_* D\rtimes E$, and we can
provide $A$ with a grading compatible with the equivalence (but we don't know
if the grading can be chosen to be non-negative).

\begin{rem}
The gradings on blocks with abelian defect should satisfy some Koszulity
properties (cf \cite{Pe}, as well as work of Chuang). Turner \cite{Tu2}
expects that gradings will even exist for blocks of symmetric
groups with non abelian defect.
\end{rem}

\begin{rem}
Using the equivalences in \S \ref{secsl2}, we obtain gradings on blocks of
abelian defect of
symmetric groups and on blocks of Hecke algebras over $\BC$. One can expect
the corresponding graded Cartan matrices to be given in terms of
Kazhdan-Lusztig polynomials. So, the equivalences carry some ``geometric
meaning''.
\end{rem}

\subsection{Dimensions}
\label{secdim}
\subsubsection{Definition and bounds}
Let us explain how to associate a dimension to a triangulated category
$\CT$ (cf \cite{dimension}). For the derived category of a finite dimensional
algebra, this is related to the Loewy length and to the global dimension, none
of which are invariant under derived equivalences.

\smallskip
Given $\CI_1$ and $\CI_2$ two subcategories of $\CT$, we denote by
$\CI_1\ast\CI_2$
the smallest full subcategory of $\CT$ closed under direct summands
and containing the objects $M$ such that there is a distinguished triangle
$$M_1\to M\to M_2\rightsquigarrow$$
with $M_i\in\CI_i$. Given $M\in\CT$, we denote by $\langle M\rangle$ the
smallest full subcategory of $\CT$ containing $M$ and
closed under direct summands, direct sums, and shifts.
Finally, we put $\langle M\rangle_0=0$ and define inductively
$\langle M\rangle_i=\langle M\rangle_{i-1}\ast \langle M\rangle$.

\smallskip
The dimension of $\CT$ is defined to be the smallest integer $d\ge 0$ such that
there is $M\in\CT$ with $\CT=\langle M\rangle_{d+1}$ (we set $\dim\CT=\infty$
if there is no such $d$). The notion of finite-dimensionality corresponds
to Bondal-Van den Bergh's property of being strongly finitely generated
\cite{BoVdB}.

\smallskip
Given a right coherent ring $A$, then
$\dim D^b(A)\le \gldim A$ (cf \cite[Proposition 2.6]{KrKu} and
\cite[Propositions 7.4 and 7.24]{dimension}).

Let $A$ be a finite dimensional algebra over a field $k$.
Denote
by $J(A)$ the Jacobson radical of $A$. The Loewy length of $A$ is the smallest
integer $d\ge 1$ such that $J(A)^d=0$. We have
$\dim D^b(A)<\ll(A)$.

\smallskip
Let $X$ be a separated scheme of finite type over a perfect field $k$.

\begin{thm}
We have $\dim D^b(X)<\infty$.
\begin{itemize}
\item
If $X$ is reduced, then $\dim D^b(X)\ge \dim X$.
\item
If $X$ is smooth and quasi-projective, then $\dim D^b(X)\le 2 \dim X$.
\item
If $X$ is smooth and affine, then $\dim D^b(X)=\dim X$.
\end{itemize}
\end{thm}

There doesn't seem to be any known example of a smooth projective variety
$X$ with $\dim D^b(X)>\dim X$, although this is expected to happen, for
example when $X$ is an elliptic curve (note nevertheless
that $\dim D^b(\BP^n)=n$).

Note that a triangulated category
with finitely many indecomposable objects up to isomorphism has dimension $0$.
This applies to $D^b(kQ)$, where $Q$ is a quiver of type ADE. This applies
also to the orbit categories constructed by Keller (cf \cite[\S 4.9]{KeICM},
\cite[\S 8.4]{Ke1}). They depend on a positive integer $d$, and they
are Calabi-Yau of dimension $d$.

\subsubsection{Representation dimension}
Auslander \cite{Au} introduced a measure for how far an algebra is from
being representation finite. The example of exterior algebras below shows that
this notion is pertinent.
Let $A$ be a finite dimensional algebra. The representation dimension of
$A$ is $\inf\{\gldim(A\oplus A^*\oplus M)\}_{M\in A\mMod}$. This is known
to be finite \cite{Iy}.

In \cite{repdim}, we show that this notion is related to the notion of
dimension for associated triangulated categories. For example,
$\dim D^b(A)\le \repdim A$. 

Let $A$ be a non semi-simple self-injective $k$-algebra. We have
$$2+\dim A\mstab\le \repdim A\le \ll(A)$$
(the second inequality comes from \cite[\S III.5, Proposition]{Au}).

The following Theorem is obtained by computing $\dim \Lambda(k^n)\mstab$ via
Koszul duality. It gives the first examples of algebras with representation
dimension $>3$.

\begin{thm}
Let $n$ be a positive integer.
We have
$\repdim \Lambda(k^n)=n+1$.
\end{thm}

\begin{rem}
One can actually show more quickly \cite{KrKu} that the algebra with quiver
$$
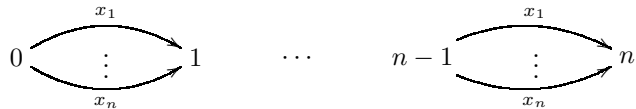
\xymatrix{
0\ar@/_1pc/[rr]_-{x_n} \ar@/^1pc/[rr]^-{x_1} & \vdots &1 & \cdots &
n-1\ar@/_1pc/[rr]_-{x_n} \ar@/^1pc/[rr]^-{x_1} & \vdots & n
}
$$
and relations $x_ix_j=x_jx_i$
has representation dimension $n$, using that its derived category is equivalent
to $D^b(\BP^n)$ \cite{Bei}.
\end{rem}

\smallskip
Using the inequality above, one obtains the following Theorem, which
solves the prime $2$  case of a conjecture of Benson.

\begin{thm}
Let $G$ be a finite group and $k$ a field of characteristic $2$.
If $G$ has a subgroup isomorphic to $(\BZ/2)^n$, then $n<
\ll(kG)$.
\end{thm}

\section{Categorifications}
\label{seccat}
This chapter discusses the categorifications of two structures, which are
related to derived equivalences. We hope these categorifications will
eventually lead to the construction of four-dimensional quantum field
theories (as advocated in \cite{CrFr}), via the construction of appropriate
tensor structures.

\subsection{$\Gsl_2$}
\label{secsl2}
\subsubsection{Abelian defect conjecture for symmetric and general linear
groups}
Let $G$
be a symmetric group and $B$ an $\ell$-block of $kG$ with defect group $D$.
Assume $D$ is abelian and let $w=\log_\ell |D|$. In 1992, a 
three steps strategy was proposed for Conjecture \ref{ADCstrong}
(inspired by the simpler character-theoretic part \cite{Routhese}):
\begin{itemize}
\item
Rickard equivalence between
$k(\BZ/\ell\rtimes \BZ/(\ell-1))\wr \GS_w$ and the principal block of
$k\GS_\ell\wr \GS_w$
\item
Morita equivalence between the principal block of
$k\GS_\ell\wr \GS_w$ and $B_w$
\item
Rickard equivalence between
$B_w$ and $B$.
\end{itemize}
Here, $B_w$ is a certain $\ell$-block of symmetric groups (a ``good block'').
Scopes \cite{Sc} has constructed a number of Morita equivalences between
blocks of symmetric groups. For fixed $w$, there are only finitely many
classes of blocks of symmetric groups up to Scopes equivalence, and
$B_w$ is defined to be
the largest block that is not Scopes equivalent to a smaller block.

The first equivalence is deduced from an equivalence between
the principal blocks of $\GS_\ell$ and $\BZ/\ell\rtimes \BZ/(\ell-1)$ via
Clifford theory \cite{Ma}.

The second equivalence was established by Chuang and Kessar \cite{ChKe},
the functor used is a direct summand of the induction functor.

The third equivalence is part of the general problem, raised by Brou\'e,
 of constructing
Rickard equivalences between two blocks of symmetric groups with
isomorphic defect groups (equivalently, with same local structure).
Rickard \cite{RiSn} constructed complexes of bimodules that he conjectured would
solve that problem, generalizing
Scopes construction (case where the complex has only one non-zero term).
Rickard proved the invertibility of
his complexes when they have two non-zero terms. The general case has
proven difficult to handle directly.

\begin{rem}
The same strategy applies for general linear groups (in non-describing
characteristic). Theorem \ref{jordan} reduces the study to unipotent blocks.
Step 2 above was handled in \cite{Mi,Tu1}.
As pointed out by H.~Miyachi, this generalizes Puig's
result \cite{PuLuminy} ($\GL_n(q)$, $\ell|(q-1)$).
\end{rem}

\begin{rem}
``Good'' blocks of symmetric groups have ``good'' properties. After the Morita
equivalence Theorem of \cite{ChKe}, their properties
were first analyzed by Miyachi \cite{Mi}, in the more complicated case
of general linear groups: decomposition matrices and radical series
of Specht modules were determined in the abelian defect case, by a direct
analysis of the wreath product. As a
consequence, decomposition matrices were known for good blocks of
Hecke algebras in characteristic zero.
For good blocks of symmetric groups with abelian defect, as well as for
Hecke algebras in characteristic zero,
a direct computation of the decomposition numbers is given in
\cite{JaLyMa} (cf also \cite{JaMa} for earlier results in that direction)
and another approach is the
determination of the relevant part of the canonical/global crystal basis
\cite{ChTa1, ChTa2, LeMi}.

For blocks of symmetric groups
with non abelian defect, the decomposition matrices can be described in
terms of decomposition matrices of smaller symmetric groups and remarkable
structural properties are conjectured by Turner \cite{Tu2,Tu3,Pa}. Good
blocks have also been used by Fayers
for the classification of irreducible Specht
modules \cite{Fa1} and to show that blocks of weight $3$ have decomposition
numbers $0$ or $1$ (for $\ell>3$) \cite{Fa2}.
\end{rem}

\subsubsection{Fock spaces}
\label{secFock}
Let us recall the Lie algebra setting for symmetric group representations
(cf e.g. \cite{Ar2}).
Let $F=\bigoplus_{n\ge 0} \BQ\otimes_\BZ K_0(\BC\GS_n\mMod)$. The complex
irreducible representations of the symmetric group $\GS_n$ are parametrized
by partitions of $n$ and we obtain a basis of $F$ parametrized by all
partitions. We view $F$ as a Fock space, with an action of $\hat{\Gsl}_\ell$
and we recall a construction of this action, for the generators
$e_a$ and $f_a$ (where $a\in\BF_\ell$).

We have a decomposition 
$$\Res^{\BF_\ell\GS_n}_{\BF_\ell\GS_{n-1}}=
\bigoplus_{a\in\BF_\ell} F_a,$$
where $F_a(M)$ is the generalized $a$-eigenspace of
$X_n=(1,n)+(2,n)+\cdots+(n-1,n)$. Taking classes in $K_0$ and summing
over all $n$, we obtain endomorphisms $f_a$ of
$$V=\bigoplus_{n\ge 0} \BQ\otimes_\BZ K_0(\BF_\ell\GS_n\mMod).$$
Using induction, we obtain similarly endomorphisms $e_a$ (adjoint to the
$f_a$).
The decomposition lifts to a decomposition of
$\Res^{\BZ_\ell\GS_n}_{\BZ_\ell\GS_{n-1}}$ and
we obtain endomorphisms $e_a$ and $f_a$ of $F$.
The decomposition map $F\to V$ and the Cartan map
$\bigoplus_{n\ge 0} \BQ\otimes_\BZ K_0(\BF_\ell\GS_n\mproj)\to F$ 
are morphisms of $\hat{\Gsl}_\ell$-modules. The image of the Cartan map is the
irreducible highest weight submodule $L$ of $F$ generated by $[\emptyset]$.

Let us note two important properties relating the module structure of
$V$ and the modular representation theory of symmetric groups:
\begin{itemize}
\item The decomposition of $V$ into weight spaces corresponds to the block
decomposition.
\item Two blocks have isomorphic defect groups if and only if they are
in the same orbit under the adjoint action of the affine Weyl group
$\tilde{A}_{\ell-1}$.
\end{itemize}

In order to prove that two blocks of symmetric groups with isomorphic
defect groups are derived equivalent, it is enough to consider a block
and its image by a simple reflection $s_a$ of $\tilde{A}_{\ell-1}$ (this
involves only the $\Gsl_2$-subalgebra generated by $e_a$ and $f_a$). This
is the situation in which Rickard constructed his complexes $\Theta_a$.

\begin{rem}
These constructions extend to Hecke algebras of symmetric groups over
$\BC$, at an $\ell$-th root of unity (here, $\ell\ge 2$ can be an arbitrary
integer).
In that situation, the classes of the indecomposable projective modules form
the canonical/global crystal
basis of $L$ (Lascoux-Leclerc-Thibon's conjecture, proven by Ariki \cite{Ar1},
cf also \cite{Gr1}).
\end{rem}

\subsubsection{$\Gsl_2$-categorifications}
We describe here joint work with
J.~Chuang \cite{ChRou} (cf also \cite{mexico} for a survey and
\cite{Gr2,GrVa,BerFreKho,FreKhoStr} for related work). This is
the special case of a more general theory under construction
for Kac-Moody algebras.

Let $k$ be an algebraically closed field and $\CA$ a $k$-linear abelian
category all of whose objects have finite composition series.

An $\Gsl_2$-categorification on $\CA$ is the data of
\begin{itemize}
\item $(E,F)$ a pair of adjoint exact functors $\CA\to\CA$
\item $X\in \End(E)$, $T\in\End(E^2)$, $q\in k^\times$, and $a\in k$
(with $a\not=0$ if $q\not=1$)
\end{itemize}
satisfying the following properties:

\begin{itemize}
\item $[E]$ and $[F]$ give rise to a locally finite representation of
$\Gsl_2$ on $K_0(\CA)$
\item for $S$ a simple object of $\CA$, $[S]$ is a weight vector
\item $F$ is isomorphic to a left adjoint of $E$
\item $(T\idun_E)\circ (\idun_E T)\circ (T\idun_E)=
(\idun_E T)\circ (T\idun_E)\circ (\idun_E T)$
\item $(T+\idun_{E^2})\circ (T-q\idun_{E^2})=0$
\item $T\circ (\idun_E X)\circ T=
\begin{cases}
q(X\idun_E) & \text{ if } q\not=1\\
X\idun_E -T & \text{ if } q=1
\end{cases}$
\item $X-a\idun_E$ is locally nilpotent.
\end{itemize}

From that data, we define two truncated powers $E^{(n,\pm)}$
(non-canonically isomorphic), using an affine Hecke algebra action on
$E^n$. Following Rickard, we construct a complex $\Theta$ 
with terms $E^{(i,-)}F^{(j,+)}$.

The following Theorem is proved by reduction to the case of
``minimal categorifications'', which are naturally associated to
simple representations of $\Gsl_2$.

\begin{thm}
\label{equivTheta}
$\Theta$ gives rise to self-equivalences of $K^b(\CA)$ and $D^b(\CA)$. This
categorifies the action of $\begin{pmatrix} 0&1\\-1&0\end{pmatrix}$
on $K_0(\CA)$.
\end{thm}

The construction of \S \ref{secFock} provides a structure of
$\Gsl_2$-categorification on $\CA=\bigoplus_{n\ge 0}\bar{\BF}_\ell\GS_n\mMod$
(for a given $a\in\BF_\ell$). From the previous Theorem, we deduce

\begin{cor}
Two blocks of symmetric groups with isomorphic defect groups are
splendidly Rickard equivalent.

Conjecture \ref{ADCstrong} holds for blocks of symmetric groups.
\end{cor}

This Corollary has a counterpart for $\GL_n(\BF_q)$ and $\ell{\not|}q$.

\begin{rem}
In general, there is a decomposition $\CA=\bigoplus_\lambda \CA_\lambda$
coming from the weight space decomposition of $K_0(\CA)$. There is a
categorification of $[e,f]=h$ in the form of
isomorphisms $EF_{|\CA_\lambda}\iso FE_{|\CA_\lambda}\oplus
\Id_{\CA_\lambda}^{\bigoplus \lambda}$ (for $\lambda\ge 0$).
\end{rem}

\begin{rem}
One can give a definition of $\Gsl_2$-categorifications for triangulated
categories and the definition above becomes a theorem that says that
there is an induced categorification on $K^b(\CA)$ (and on $D^b(\CA)$).
\end{rem}

\begin{rem}
One can also construct $\Gsl_2$-categorifications on category $\CO$ for
$\Ggl_n(\BC)$ and for rational representations of $\GL_n(\bar{\BF}_p)$. One
deduces from Theorem \ref{equivTheta} that blocks with the same
stabilizers
under the affine Weyl groups are derived equivalent (a conjecture of
Rickard).
\end{rem}

\begin{rem}
The endomorphism $X$ has different incarnations: Jucys-Murphy element,
Casimir,...
\end{rem}

\begin{rem}
It is expected that the functors $\Theta_a$ constructed for $a\in\BF_\ell$
provide an action of the affine braid
group $B_{\tilde{A}_{\ell-1}}$ on $\bigoplus_n D^b(\BF_\ell\GS_n)$.
\end{rem}

\subsection{Braid groups}
\label{secbraid}
\subsubsection{Definition}
We present here a categorification of braid groups associated to Coxeter
groups, following \cite{mexico}.
This should be useful for the study of categories of representations
of semi-simple Lie algebras, affine Lie algebras, simple algebraic groups
over an algebraically closed field,... On the other hand, work of Khovanov
\cite{Kh} shows its relevance for invariants of links (type $A$).

Let $(W,S)$ be a Coxeter group, with $S$ finite. Let $V$ be its reflection
representation over $\BC$ and let $B_W$ be the braid group of $W$.
Let $A=\BC[V]$. Given $s\in S$, let
$F_s=0\to A\otimes_{A^s}A\xrightarrow{\text{mult}}A\to 0$,
where $A$ is in degree $1$.
This is an invertible object of $K^b(A\otimes A)$. Given two decompositions
of an element of $B_W$ in a product of the generators and their inverses,
we construct a canonical isomorphism between the corresponding products
of $F_s$. The system of isomorphism coming from the various decompositions
of an element $b\in B_W$ is transitive and, taking its limit, we obtain
an element $F_b\in K^b(A\otimes A)$. The full subcategory of
$K^b(A\otimes A)$ with objects the $F_b$'s defines a strict monoidal
category $\CB_W$.

We expect that there is a simple presentation of $\CB_W$ by generator
and relations (or rather of a related $2$-category involving subsets of $S$).
This should be related to the vanishing of certain $\Hom$-spaces, for example
$\Hom_{K^b(A\otimes A)}(F_b,F_{b'}^{-1}[i])$ should be $0$ when
$b$ and $b'$ are the canonical lifts of distinct elements of $W$.

\begin{rem}
The bimodules obtained by tensoring the $A\otimes_{A^s}A$ are Soergel's
bimodules. Soergel showed they categorify the Hecke algebra of $W$. He also
conjectured that the indecomposable objects correspond to the Kazhdan-Lusztig
basis of $W$ \cite{Soe1,Soe2}.
\end{rem}

\begin{rem}
When $W$ is finite,
one can expect that there is a construction of $\CB_W$ that does not
depend on the choice of $S$. Such a construction might then
make sense for complex reflection groups.
\end{rem}

\subsubsection{Representations and geometry}
Let $\Gg$ be a complex semi-simple Lie algebra with Weyl group $W$ and
let $\CO_0$ be the principal block of its category $\CO$. It has been
widely noticed that there is a weak action of $B_W$ on $D^b(\CO)$, using
wall-crossing functors. We show that there is a genuine action of $B_W$ on
$D^b(\CO_0)$ and there is a much more precise statement: there is a monoidal
functor from $\CB_W$ to the category of self-equivalences of $D^b(\CO_0)$.
This has a counterpart for the derived category of $B$-equivariant sheaves
on the flag variety (in which case the genuine action of the braid group
goes back to \cite{De}). These actions are compatible with
Beilinson-Bernstein's equivalence. Conversely, a suitable presentation
of $\CB_W$ by generators and relations should provide a quick proof of
that equivalence (and of affine counterparts), in the spirit of Soergel's
construction. The representation-theoretic and the geometrical categories
should be viewed as two realizations of the same ``$2$-representation''
of $\CB_W$. Also, this approach should give a new proof of the results
of \cite{AnJaSoe} comparing quantum groups at roots of unity and 
algebraic groups in characteristic $p$.

\end{document}